\documentclass[12pt]{article}
\usepackage[numbers,sort&compress]{natbib}
\usepackage{amssymb,amsmath,mathrsfs,amsfonts,amsthm}
\usepackage{enumerate}
\usepackage{paralist}
\usepackage{color}
\usepackage{epsfig}
\usepackage{graphicx}
\usepackage{setspace}
\usepackage{appendix}
\begin{document}
 \def\pd#1#2{\frac{\partial#1}{\partial#2}}
\def\dfrac{\displaystyle\frac}
\let\oldsection\section
\renewcommand\section{\setcounter{equation}{0}\oldsection}
\renewcommand\thesection{\arabic{section}}
\renewcommand\theequation{\thesection.\arabic{equation}}

\newtheorem{thm}{Theorem}[section]
\newtheorem{cor}[thm]{Corollary}
\newtheorem{lem}[thm]{Lemma}
\newtheorem{prop}[thm]{Proposition}
\newtheorem*{con}{Conjucture}
\newtheorem*{questionA}{Question}
\newtheorem*{thmA}{Theorem A}
\newtheorem*{thmB}{Theorem B}
\newtheorem{remark}{Remark}[section]
\newtheorem{definition}{Definition}[section]

\title{Optimisation of total population  in logistic model with nonlocal dispersals  and heterogeneous environments
\thanks{The second author is  supported by  NSF  of China (No. 11971498). } }

\author{Xueli Bai\\ {School of Mathematics and Statistics, Northwestern Polytechnical University,}\\{\small 127 Youyi Road(West), Beilin 710072,
Xi'an, P. R. China.}\\
Fang Li{\thanks{Corresponding author.
E-mail: lifang55@mail.sysu.edu.cn}}\\ {School of Mathematics, Sun Yat-sen
University,}\\{\small No. 135, Xingang Xi Road, Guangzhou 510275, P. R. China.} \\
Maolin Zhou \\ {Chern Institute of Mathematics and LPMC, Nankai University,}\\{\small Tianjin 300071,  P. R. China.}  }

\date{}
\maketitle{}

\begin{abstract}
In this paper, we   investigate the issue of maximizing the total equilibrium population  with respect to  resources distribution $m(x)$ and diffusion rates $d$ under the prescribed total amount of resources in a logistic model with nonlocal dispersals.
Among other things,  we show that for $d\geq 1$,  there exist $C_0, C_1>0$, depending on $\|m\|_{L^1}$ only,  such that
$$
C_0\sqrt{d} \leq \textrm{  supremum of    total population  }\leq C_1 \sqrt{d}.
$$
However,  when replaced by random diffusion,  a conjecture, proposed by Ni and justified in \cite{BHL}, indicates that in the one-dimensional case,
$$
\textrm{  supremum of    total population  } = 3\|m\|_{L^1}.
$$
This reflects  serious discrepancies between models with local and nonlocal dispersal strategies.
\end{abstract}



{\bf Keywords} total population, nonlocal dispersal, heterogeneity
\vskip3mm {\bf MSC (2010)}: 35K57, 92D25, 45K05

\section{Introduction}
Dispersal is an important feature of  life histories of many organisms and often crucial for their persistence. 
Understanding the effect of dispersal in heterogeneous environment on population
dynamics is an important issue in spatial ecology \cite{CC-book}.
Total population is an important indicator for  persistence of species. If the quantity is at low level, the risk of extinction will increase, while if the quantity is at high level, it will lead to  shortage of resources and intense pressure of competition, which may japodize the existing stability of the multi-species systems \cite{LouCN}.
Therefore, {\it an interesting problem in spatial ecology is how     dispersal strategies of the species and the distribution of resources affect the total population.}

Our study  is motivated by a series of intriguing questions and work related to total equilibrium population in  a single logistic equation with random diffusion as follows
\begin{equation}\label{general single}
\begin{cases}
u_t =   d\Delta u    +   u    [m(x)-  u]  & x\in \Omega,\ t>0,\\
\frac{\partial u}{\partial \nu} =0  & x\in \partial\Omega,\ t>0,
\end{cases}
\end{equation}
where  $u$ represents the population density of a species at location $x\in\Omega$ and at time $t>0$, $d$ is the dispersal rate of the species which is assumed to be a positive constant, the habitat $\Omega$ is a bounded domain in $\mathbb R^n$ and $\nu$ denotes the unit outward normal vector.
The function $m(x)$ is the intrinsic growth rate or carrying capacity, which reflects the environmental influence on the species $u$. Unless designated otherwise, we assume that $m(x)$ satisfies the following condition:
\medskip\\
\noindent\textbf{(M)} \hspace{0.5cm}  $ m(x)\in L^{\infty}(\Omega),\ m(x)\geq 0 \ \textrm{and}\  m\not\equiv \textrm{const} \ \textrm{on} \  \bar\Omega.$
\medskip\\
It is known that if $m$ satisfies the assumption \textbf{(M)}, then for every $d>0$, the problem (\ref{general single}) admits a unique positive steady state, denoted by $\theta_{d,m}(x)$, which is globally asymptotically stable (see e.g. \cite{CC-book}). In addition, a remarkable property concerning $\theta_{d,m}(x)$ was first observed in \cite{Lou06}
\begin{equation}\label{observation}
\int_{\Omega}\theta_{d,m}(x)\,dx > \int_{\Omega}m(x)\,dx \quad\text{ for all }d > 0.
\end{equation}
Biologically, this  indicates that when coupled with diffusion, a heterogeneous environment can support a total population larger than the  total carrying capacity of the environment, which is quite different from homogeneous environment.  Simply speaking, heterogeneity  of resources can benefit survival of species. This theory is  further  confirmed experimentally \cite{ZDN}.
Moreover, it is known that (see e.g. \cite{Lou06})
\begin{equation}\label{observation2}
\lim_{d\rightarrow 0^+} \int_{\Omega} \theta_{d,m}(x)\, dx = \lim_{d\rightarrow \infty} \int_{\Omega} \theta_{d,m}(x)\, dx = \int_{\Omega}m(x)\,dx.
\end{equation}
This indicates that for given $m(x)$, the  total   population as a function of the diffusion rate $d$ is not monotone and achieves its maximum at some intermediate value. Examples constructed in \cite{LiangLou12} show that the local maximum  might not be unique.
These observations  naturally lead to a biological question:

\medskip
\noindent\textbf{Question.}  Given the total amount of resources, how should we distribute the resource and/or adjust diffusion rate to maximize the total equilibrium population?
\medskip

For simplicity, denote
$$
\mathcal{M}_1=\{ m \,|\, m \text{ satisfies condition \textup{(\textbf{M})}},\ \int_{\Omega}m(x) dx =1\}.
$$
In the one-dimensional case,  W.-M. Ni
conjectured that the supremum of the total population over all $d>0$ and $m\in \mathcal{M}_1$ is $3$. This conjecture
is confirmed in \cite{BHL}. However, for higher dimensional case,  it is proved in \cite{Jumpei} that the supremum is unbounded.
Moreover,  when the  diffusion rate is fixed, optimisation of total population size with respect to  resources distributions is studied and existing results indicate that the optimal configuration is of bang-bang type.  See \cite{NY2018, MNP, Mazari2021, MNP2022}.  This question is also studied in patchy environment \cite{NLouY2021}.

The main purpose of this paper  is to investigate this question  for the following single species model with nonlocal dispersal strategy
\begin{equation}\label{main single}
\begin{cases}
  u_t(x,t)=d\mathcal {L} [u]     +  u   [m(x)-  u]   &x\in \Omega,\ t>0,\\
  u(x,0)=u_0\geq 0, &x\in \Omega,
\end{cases}
\end{equation}
where the nonlocal diffusion operator $\mathcal{L}$, which corresponds to nonlocal homogeneous Neumann boundary condition, is defined as:
 \begin{equation}\label{nonlocal operator}
   \mathcal{L} [u] :=\int_\Omega k(x,y)u(y)dy-\int_{\Omega} k(y,x)dy u(x),
 \end{equation}
and the dispersal kernel function $k(x,y)\geq 0$ describes the probability to jump from one
location to another. 
Although the most popular forms of continuum
models have been those given by differential equations (local diffusion), in many situations in ecology (e.g. \cite{Cain,Clark1,Clark2,schurr}), dispersal is better described as a long range process rather than  a local one, and integral operators appear as a natural choice. Depending on the background, there are many forms of nonlocal models. For a detailed introduction of modeling, see the book \cite{Du-book} and the survey paper \cite{DuICM}.  The nonlocal diffusion operator studied in this paper  is a  general form which appears commonly in different types of models in ecology. See \cite{Allen1996, HMMV, Kot1996, Lee2001, Lutscher, Medlock2003, Meysman2003, Mogilner-E, Othmer} and the references therein.



From now on,  assume that the kernel $k$ satisfies
\medskip\\
\noindent\textbf{(K)} \   $k(x,y)\in C(\mathbb R^n\times \mathbb R^n)$ is nonnegative and  $k(x,x)>0$ in $\mathbb R^n$. $k(x,y)$ is symmetric,

\ \ \ \ i.e., $k(x,y)=k(y,x)$. Moreover,  $\int_{\mathbb R^n} k(x,y)dy =1$.

\medskip
\noindent and for simplicity, denote
$$
a(x) = \int_{\Omega} k(y,x) dy\leq 1.
$$

First of all, we prepare the existence and uniqueness result for the model (\ref{main single}) provided that $m(x) \in L^{\infty} (\Omega)$.

\begin{thm}\label{thm-existence}
Assume that $m(x) \in L^{\infty} (\Omega)$ is nonconstant and the kernel $k$ satisfies \textbf{(K)}. Define
\begin{equation}\label{mu0}
\displaystyle\mu_0=\mu_0(m) =\sup_{0\neq\psi\in L^2(\Omega)} \frac{\int_{\Omega}\left(d\mathcal {L} [\psi] (x)\psi(x) +m(x)\psi^2(x)\right)dx}{\int_{\Omega} \psi^2(x) dx}.
\end{equation}
Then  the problem (\ref{main single})
admits a unique positive steady state in $L^{\infty}(\Omega)$  if and only if $\mu_0>0$. In particular, if $m\in\mathcal{M}_1$, then the problem (\ref{main single})
admits a unique positive steady state, denoted by $\theta_{d,m}$, in $L^{\infty}(\Omega)$  for any $d>0$.
\end{thm}

When $m\in C(\bar\Omega)$, the existence and uniqueness of positive steady state for the model (\ref{main single}) has been studied thoroughly.   See \cite{BZh} for symmetric operators in the one dimensional case and  \cite{BL2, Coville2010}  for nonsymmetric operators. The proofs of these studies rely on the properties of nonlocal eigenvalue problems, thus the condition $m\in C(\bar\Omega)$ is required. However, to study the questions in this paper, the condition $m(x) \in L^{\infty} (\Omega)$ is necessary.
For this purpose, we develop a different approach, which depends on the application of energy functional.
It is known that the properties (\ref{observation}) and (\ref{observation2}) also hold for the solution  to the nonlocal problem (\ref{main single}) provided that  $m$ satisfies the assumption \textbf{(M)} (see e.g. \cite{SuLi}).

Our first result indicates that   the supremum of total equilibrium population $\displaystyle \int_{\Omega}\theta_{d,m}\,dx$   over $m\in\mathcal{M}_1$ and $d\geq 1$ is of order $\sqrt{d}$, where the existence of $\theta_{d,m}$ is guaranteed by Theorem \ref{thm-existence}.

\begin{thm}\label{theorem-unbounded-rate}
	Assume that $\Omega$ is a bounded domain in $\mathbb R^n$, $n\geq 1$. Then   there exist $C_0,\, C_1>0$, independent of $m\in\mathcal{M}_1$,  such that for $d\geq 1$
	\begin{equation}\label{thm-upper-lower}
		C_0 \sqrt{d} \leq \sup\,\left\{ \int_{\Omega}\theta_{d,m}\,dx \, \Big|\, m\in\mathcal{M}_1 \right\} \leq  C_1 \sqrt{d},
	\end{equation}
where $\theta_{d,m}$ denotes the unique positive steady state to the problem (\ref{main single}).
\end{thm}

Recall that the studies in \cite{BHL, Jumpei} for the local model (\ref{general single}) reveal that  the supremum of the total equilibrium population over all $d>0$ and $m\in \mathcal{M}_1$ is $3$  for the one-dimensional case and unbounded  for the  higher dimensional case, respectively.  However, Theorem \ref{theorem-unbounded-rate}   indicates that the supremum is always unboundedness in any dimensional case due to the unboundedness of diffusion rate $d$ and $\|m\|_{L^{\infty}}$.  This demonstrates  serious discrepancies between models with local and nonlocal dispersal strategies.

On the basis of Theorems \ref{theorem-unbounded-rate}, we further explore how to combine  resources distribution $m(x)$ and diffusion rate $d$ such that
the total equilibrium population of the problem (\ref{main single}) is of order $\sqrt{d}$ as $d$ goes to infinity.
The following result provides an equivalent characterization.

\begin{thm}\label{thm-optimal}
Assume that $\Omega$ is a bounded domain in $\mathbb R^n$, $n\geq 1$, $\theta_{d,m}$ denotes the unique positive steady state to the problem (\ref{main single}) and $\mathcal{S}$ is a subset of $ \{d>0 \}\times \mathcal{M}_1$. Then there exists $C>0$ such that $\displaystyle \int_{\Omega}\theta_{d,m}\,dx \geq  C \sqrt{d}$  for all $(d,m)\in\mathcal{S}$ if and only if the following assumption is valid in $\mathcal{S}$:
\medskip\\
\noindent\textbf{(A)} \ \   There exist $\varepsilon_0,\, D>0$ such that for any $(d,m)\in\mathcal{S}$ with $d>D$,  $m$ satisfies
$$
\displaystyle\int_{\{{m\over d}>(1+\varepsilon_0)a\}} m(x) dx \geq \varepsilon_0,
$$
where
$$
a(x) = \int_{\Omega} k(y,x) dy\leq 1.
$$
\end{thm}

\medskip

Intuitively, Theorem \ref{thm-optimal} demonstrates that given  total amount of  resources,  the total equilibrium population could reach the order $\sqrt{d}$ as $d\rightarrow \infty$ if and only if certain amount  of resources concentrates at certain height related to diffusion rate as described in the assumption \textbf{(A)}.
As mentioned earlier, a series of results for the local model (\ref{general single}) in  \cite{NY2018, MNP, Mazari2021, MNP2022} indicate that  when the diffusion rate is fixed, bang-bang type is the optimal configuration of the resources for maximizing total population size.
Therefore, generally speaking, for both local and nonlocal models, concentration of resources is benefit for survival of species.

We emphasize that Theorem \ref{thm-optimal}  provides an equivalent   criterion to determine whether the total equilibrium population could reach the order $\sqrt{d}$ as $d\rightarrow \infty$ under the prescribed total amount of  resources.
With the help of  this  equivalent   criterion,  some concrete examples are constructed
in Section 4.2 to elaborate how to choose  $m$  for $d$ large such that the corresponding total population could reach the optimal order $\sqrt d$ and  how properties of the kernel functions and the locations where the resources concentrate affect the total population.


A surprising and important feature in Theorems \ref{theorem-unbounded-rate} and \ref{thm-optimal} is that the estimates for  the total population are independent of $\|m\|_{L^{\infty}}$, although according to  the estimate $\|\theta_{d,m}\|_{L^{\infty}}\le\|m\|_{L^{\infty}}$,   the unboundedness of $\|m\|_{L^{\infty}}$ is necessary for unboundedness of the total population. Indeed, we only require that the total resources, i.e., $\|m\|_{L^1}$, is fixed. Biologically, this is the most natural and basic assumption. Mathematically, it is quite challenging to obtain these estimates, since we could only rely on the $L^1$ norm of $m$ and need avoid the appearance of other norms.
Throughout the proofs of Theorems \ref{theorem-unbounded-rate} and \ref{thm-optimal},  the key idea is to divide the domain $\Omega$ in a proper way such that we could establish more precise estimates about $\theta_{d,m}$, $m$ in each sub-domain.

At the end, notice that for Theorems \ref{thm-existence}--\ref{thm-optimal}, only the assumption \textbf{(K)} is imposed on the kernel function $k(x,y)$.  However, if the kernel function $k(x,y)$ is of the special form
$$k(x,y)=J(x-y),$$
it is proved that solutions of nonlocal equations with suitably rescaled  kernel functions converge to solutions of local equations \cite{Rossi}.
In particular, if $J$ also satisfies
$$
\int_{\mathbb R^n} |x|^2J(x)dx<+\infty.
$$
then the convergence relation between solutions to equations with  rescaled nonlocal operator and random diffusion is verified \cite{Rossi}.
This reflects the complexity and diversity of nonlocal operators considered in this paper.
Therefore, for nonlocal models, it is possible that bang-bang type is not optimal location of resources maximizing total population  for fixed diffusion rate.
The answer might depend on the specific properties  of kernel functions $k(x,y)$.
Additionally, Theorems \ref{theorem-unbounded-rate}  and \ref{thm-optimal} are related to the combined effects of diffusion rate and resources  distribution  on maximizing total equilibrium population. This remains unknown for local models. We will return to these problems in future work.



This paper is organized as follows. Theorem \ref{thm-existence} is   proved in Section 2. Section 3 is devoted to the proof of  Theorem  \ref{theorem-unbounded-rate}.  At the end, the proof of  Theorem \ref{thm-optimal} is present and some concrete examples are discussed in Section 4.

\section{Existence and uniqueness of positive steady state}
In this section, we establish  the existence and uniqueness of positive steady state to the  problem of (\ref{main single}) when $m(x) \in L^{\infty} (\Omega)$.

\begin{proof}[Proof of Theorem \ref{thm-existence}]
First, if   the problem (\ref{main single}) admits a positive steady state, denoted by $\theta$, in $L^{\infty}(\Omega)$, then it is easy to see that $\mu_0>0$ by choosing $\psi=\theta$.

The rest of the proof is devoted to proving the other direction.
Let $u$ be the solution of
\begin{equation}\label{single-M}
\begin{cases}
  u_t=d\mathcal {L} [u] (x,t)    +  u(x,t)   [m(x)-  u(x,t)]   &x\in \Omega,\ t>0,\\
  u(x,0)=\|m\|_{L^{\infty}} &x\in \Omega.
\end{cases}
\end{equation}
 Thus, $u$ is decreasing in $t$ and there exists $\theta^*\in L^\infty(\Omega)$ such that $u(x,t)\rightarrow \theta^*(x)$ pointwisely as $t\rightarrow\infty$.  Moreover, $\theta^*$ is a steady state of (\ref{single-M}).

Now we show that $\theta^*\not\equiv 0$. Suppose that it is not true, that is $u(x,t)\rightarrow 0$ pointwisely  as $t\rightarrow\infty$.
Since $\mu_0>0$, by the definition of $\mu_0$ we can choose $0\not\equiv\psi_0\in L^2$ such that
\begin{equation}\label{eigen2}
\int_{\Omega}\left(d\mathcal {L} [\psi_0] (x)\psi_0(x) +m(x)\psi_0^2(x)\right)dx\ge \frac{\mu_0}{2}\int_{\Omega} \psi_0^2 dx>0.
\end{equation}
Let $\psi_i:=\min\{\psi_0, i\}$, obviously $\psi_i\rightarrow\psi_0$ in $L^2(\Omega)$ as $i\rightarrow\infty$. Combined with (\ref{eigen2}), we can fix $i=i_0$ large enough, such that
$$
\int_{\Omega}\left(d\mathcal {L} [\psi_{i_0}] (x)\psi_{i_0}(x) +m(x)\psi_{i_0}^2(x)\right)dx\ge \frac{\mu_0}{4} \int_{\Omega} \psi_{i_0}^2 dx > 0.
$$
Set $\phi:=\varepsilon_{i_0}\psi_{i_0}$, with $\displaystyle\varepsilon_{i_0}= \frac{1}{i_0} \min\{\|m\|_{L^{\infty}},\frac{\mu_0}{8}\}$. It is routine to verify that
\begin{equation}\label{eigne-ep}
\int_{\Omega}\left(d\mathcal {L} [\phi] (x)\phi(x) +m(x)\phi^2(x)\right)dx-\frac{2}{3} \int_{\Omega} \phi^3dx \ge [\frac{\mu_0}{4}-\varepsilon_{i_0} i_0]\int_{\Omega} \phi^2 dx> 0.
\end{equation}
Suppose that $v$ is the  solution of
$$
\begin{cases}
  v_t=d\mathcal {L} [v] (x,t) +(m-v)v   &x\in \Omega,t>0,\\
  v(x,0)=\phi   &x\in \Omega,
\end{cases}
$$
and define
$$
E[v](t):=\frac{1}{2}\int_{\Omega}\left( d\mathcal {L} [v] v + mv^2 \right) dx - \frac{1}{3} \int_{\Omega}  v^3 dx.
$$
By comparison principle $\phi\le \|m\|_{L^{\infty}}$ implies that $v\le u$. Thus, $v\rightarrow 0$ in pointwisely as $t\rightarrow\infty$, and furthermore
\begin{equation}\label{E to 0}
E[v](t)\rightarrow 0 ~{\rm as}~ t\rightarrow\infty.
\end{equation}
However, since $k(x,y)$ is symmetric, straightforward computation yields that
$$
  \frac{d}{dt}E[v](t) =\int_{\Omega} v_t^2 dx\ge 0.
$$
Together with (\ref{eigne-ep}),  one sees that $ E[v](t)$ is a increasing function with positive initial data, which contradicts to (\ref{E to 0}).

Hence $\theta^*(x)\ge 0$ is a nontrivial steady state of (\ref{main single}).
Furthermore, denote $A:=\{x\in\Omega\, |\, \theta^*(x)=0\}$.  Due to the assumption \textbf{(K)}, a contradiction can be derived easily by integrating both sides of the equation satisfied by $\theta^*$ in $A$ if $A$ has positive measure. This yields that $\theta^*>0$ a.e. in $\Omega$.

It remains to show the uniqueness  of positive steady state to the  problem of (\ref{main single}) in $L^{\infty}(\Omega)$.  Suppose that $\theta\in L^\infty(\Omega)$ is a positive steady state of (\ref{main single}), i.e. $\theta$ satisfies
$$
d\mathcal {L} [\theta] (x)    +  \theta(x)   [m(x)-  \theta(x)]=0,
$$
By multiplying both sides by $\theta^{p-1}$ and integrating over $\Omega$,
we have
\begin{eqnarray*}
&& \int_{\Omega} \theta^{p+1}(x)dx - \int_{\Omega} m(x) \theta^p(x) dx\\
&=& d \int_{\Omega} \theta^{p-1}(x) \left[\int_\Omega k(x,y)\theta(y)dy-a(x)\theta(x) \right]dx\\
&\leq &d \int_{\Omega} \theta^{p-1}(x) \left(\int_\Omega k(x,y)dy\right)^{p-1\over p}\left(\int_{\Omega}k(x,y)\theta^{p}(y)dy\right)^{1\over p}dx-d\int_{\Omega}a(x)\theta^p(x) dx\\
&\leq & d \left(\int_{\Omega} \int_\Omega k(x,y)dy\theta^p(x) dx\right)^{p-1\over p}\left(\int_{\Omega}\int_{\Omega}k(x,y)\theta^{p}(y)dy dx\right)^{1\over p}-d\int_{\Omega}a(x)\theta^p(x) dx\\
&=&d \int_{\Omega} \int_\Omega k(x,y)dy\theta^p(x) dx -d\int_{\Omega}a(x)\theta^p(x) dx\leq 0,
\end{eqnarray*}
since $k(x,y)$ satisfies the assumption \textbf{(K)} and $\displaystyle a(x) = \int_{\Omega} k(y,x)dy\leq 1$. Thus it is easy to see that
$$\|\theta\|_{L^{p+1}}\le \|m\|_{L^{p+1}},$$
which yields that
\begin{equation}\label{u-m-bound}
\|\theta\|_{L^{\infty}}\le \|m\|_{L^{\infty}},
\end{equation}
since $p$ is arbitrary.
Then thanks to (\ref{single-M}), it follows that $\theta(x)\le \theta^*(x)$.
Straightforward computation gives
\begin{align*}
  \int_{\Omega} (\theta^*-\theta)\theta\theta^*dx&=\int_{\Omega}  (m-\theta)\theta\theta^*dx-\int_{\Omega} (m-\theta^*)\theta^*\theta dx\\
  &=-d\int_{\Omega} \mathcal{L}[\theta]\theta^*dx+d\int_{\Omega} \mathcal{L}[\theta^*]\theta dx =0,
\end{align*}
which implies that $\theta\equiv \theta^*$.

Obviously, if $m\in\mathcal{M}_1$, then $\mu_0(m)>0$.
The proof is complete.
\end{proof}

\begin{remark}
When the nonlocal diffusion operator $\mathcal{L} [u] $ in the problem (\ref{main single}) is replaced by
$$
\int_\Omega k(x,y)u(y)dy-  u(x),
$$
which corresponds to nonlocal homogeneous Dirichlet boundary condition, the existence and uniqueness of positive steady state in $L^{
\infty}(\Omega)$ is also equivalent to $\mu_0>0$. The proof is
 the same as that of Theorem \ref{thm-existence}.
\end{remark}



\section{Estimates for  the supremum  of total  population}
This section is devoted to the proof of Theorems  \ref{theorem-unbounded-rate},  where
we estimate the supremum of the ratio between total equilibrium population and  total carrying capacity  over  $d>0$ and $m\in L^{\infty}$ for the nonlocal model (\ref{main single}).

The proof of Theorems  \ref{theorem-unbounded-rate} consists of two parts: the upper bound and lower bound for
$$
\displaystyle \sup\,\left\{ \int_{\Omega}\theta_{d,m}\,dx \, \Big|\, m\in\mathcal{M}_1 \right\}
$$
and for clarity, we present the proofs of these two parts in two subsections respectively.
For convenience of readers, some explanations are present as follows.
\begin{itemize}
\item In the proof of the upper bound, we divide the domain $\Omega$ based multiples of $d$ as follows
    $$
    \Omega_i = \{ x\in \Omega \ | \ \theta_{d,m}>K_i d  \}, \ i=1,2,
    $$
    where the constants $K_i$, i=1,2, will be chosen carefully such that we first estimate
    $$
    \int_{\Omega_1} \theta_{d,m}(x) dx,\ \  \int_{\Omega\setminus\Omega_1} \theta_{d,m}(x) dx
    $$
    separately to obtain an upper bound of the total population in the order of $d$,  and then improve the upper bound to the order of $\sqrt d$ by estimating
    $$
    \int_{\Omega_2} \theta_{d,m}(x) dx,\ \  \int_{\Omega\setminus\Omega_2} \theta_{d,m}(x) dx
    $$
    separately.
\item In the proof of the lower bound, we construct examples to demonstrate that the order $O(\sqrt{d})$ could be achieved for $d$ large under the prescribed total carrying capacity.
\end{itemize}

\subsection{Upper bound for $\displaystyle \sup\,\left\{ \int_{\Omega}\theta_{d,m}\,dx \, \Big|\, m\in\mathcal{M}_1 \right\} $}
In this subsection, we show that there exists $C_1$, independent of $ m\in\mathcal{M}_1$,  such that
\begin{equation}\label{thm-upper}
\sup\,\left\{ \int_{\Omega}\theta_{d,m}\,dx \, \Big|\, m\in\mathcal{M}_1 \right\} \leq  C_1 \sqrt{d}.
\end{equation}

\begin{proof}[Proof of the upper bound (\ref{thm-upper})]

Thanks to Theorem \ref{thm-existence}, when $m(x)$ satisfies the condition \textbf{(M)}, the problem (\ref{main single}) always admits a unique positive steady state, denoted by $\theta_{d,m}$, i.e., $\theta_{d,m}$ satisfies
\begin{equation}\label{single-ss}
d\left( \int_{\Omega} k(x,y) \theta(y) dy - a(x) \theta (x)\right)   +  \theta(x)   [m(x)-  \theta(x)] =0\ \ \ x\in\Omega,
\end{equation}
where
$$
a(x) = \int_{\Omega} k(y,x) dy\leq 1.
$$

In the following proof, we  keep  $ \displaystyle\int_{\Omega}m(x)\,dx$ in the estimates to emphasize the role played by total carrying capacity, though indeed $ \displaystyle\int_{\Omega}m(x)\,dx=1$ since $m\in\mathcal{M}_1$.

.

First, we establish a rough estimate for $\displaystyle \int_{\Omega}\theta_{d,m}(x) dx $.
Set
\begin{eqnarray*}
\Omega_1 = \{x\in \Omega \ | \ \theta_{d,m}(x)>K_1 d \},\  \textrm{where}\ K_1= 2 \|k \|_{L^{\infty}} |\Omega|.
\end{eqnarray*}
For any $x\in \Omega_1$,
\begin{eqnarray}\label{pf-thm-Omega1-x}
\theta_{d,m}(x) &=& m(x) + \frac{d\left( \int_{\Omega} k(x,y) \theta_{d,m}(y) dy - a(x) \theta_{d,m} (x)\right)}{\theta_{d,m}(x)}\cr
&\leq  &  m(x) + \frac{d \|k \|_{L^{\infty}} \int_{\Omega}  \theta_{d,m}(x) dx }{K_1 d}\cr
&= &  m(x) + \frac{1 }{2 |\Omega|}\int_{\Omega}  \theta_{d,m}(x) dx,
\end{eqnarray}
which implies that
$$
\int_{\Omega_1} \theta_{d,m}(x) dx\leq  \int_{\Omega} m(x) dx + \frac{|\Omega_1| }{2|\Omega| }\int_{\Omega}  \theta_{d,m}(x) dx \leq \int_{\Omega} m(x) dx + \frac{1 }{2 }\int_{\Omega}  \theta_{d,m}(x) dx.
$$
Then for any $d>0$,
\begin{equation*}
\int_{\Omega} \theta_{d,m}(x) dx = \int_{\Omega_1} \theta_{d,m}(x) dx+ \int_{\Omega\setminus\Omega_1} \theta_{d,m}(x) dx\leq \int_{\Omega} m(x) dx + \frac{1 }{2 }\int_{\Omega}  \theta_{d,m}(x) dx + K_1 d|\Omega\setminus\Omega_1|.
\end{equation*}
Thus for any $d>0$,
\begin{equation}\label{remark-bounded-d}
\int_{\Omega} \theta_{d,m}(x) dx  \leq 2 \int_{\Omega} m(x) dx+ 2 K_1 d|\Omega\setminus\Omega_1|,
\end{equation}
and thus 
and for all $d\geq 1$,
\begin{equation}\label{pf-thm-roughbound}
\int_{\Omega} \theta_{d,m}(x) dx \leq 2\left(\int_{\Omega} m(x) dx +K_1 |\Omega|\right)d.
\end{equation}

Next, set
\begin{eqnarray*}
\Omega_2 =  \{x\in \Omega \ | \ \theta_{d,m}(x)>K_2 d \},\  \textrm{where}\ K_2= \frac{4\left(\int_{\Omega} m(x) dx+ K_1 |\Omega|\right)\|k \|_{L^{\infty}}}{\min_{\bar\Omega} a(x)} + 2 \|k \|_{L^{\infty}} |\Omega|
\end{eqnarray*}
and we prepare an estimate for  $|\Omega_2|$ in term of $d$. Denote
$$
\tilde\Omega_2 =\left \{ x\in \Omega \ \Big|\ m(x) \geq {d\over 2} a(x) \right\}.
$$
Obviously,
$$
\int_{\Omega} m(x) dx \geq \int_{\tilde\Omega_2} m(x) dx \geq {d\over 2} \min_{\bar\Omega} a(x)  |\tilde\Omega_2|,
$$
which implies that
$$
|\tilde\Omega_2| \leq  {1\over d} {2\over \min_{\bar\Omega} a(x) }  \int_{\Omega} m(x) dx.
$$
We claim that $\Omega_2 \subseteq \tilde\Omega_2$. If the claim is true, then one has
\begin{equation}\label{pf-thm-Omega2}
|\Omega_2| \leq {1\over d} {2\over \min_{\bar\Omega} a(x) }  \int_{\Omega} m(x) dx.
\end{equation}
To prove this claim, fix any $x\in \Omega \setminus\tilde\Omega_2$, i.e., $\displaystyle m(x) <{d\over 2} a(x)$.  Based on the equation (\ref{single-ss}),
\begin{eqnarray*}
\theta_{d,m} (x) &=&
 {1\over 2} \left[m(x) -d a(x) + \sqrt{(m(x) -d a(x))^2 + 4 d \int_{\Omega} k(x,y) \theta_{d,m}(y) dy} \right]\\
&=&\frac{  2d \int_{\Omega} k(x,y) \theta_{d,m}(y) dy}{  -m(x) +d a(x) + \sqrt{(m(x) -d a(x))^2 + 4 d \int_{\Omega} k(x,y) \theta_{d,m}(y) dy}}\\
&\leq & \frac{  2d \int_{\Omega} k(x,y) \theta_{d,m}(y) dy}{\displaystyle d a(x)}\leq \frac{  2 \|k \|_{L^{\infty}} \int_{\Omega}  \theta_{d,m}(x) dx}{  \min_{\bar\Omega} a(x)}\\
&\leq & \frac{  4\left(\int_{\Omega} m(x) dx+ K_1 |\Omega|\right)\|k \|_{L^{\infty}}}{\min_{\bar\Omega} a(x)} d,
\end{eqnarray*}
where the last inequality is due to (\ref{pf-thm-roughbound}). Hence $\theta_{d,m}(x) < K_2 d$, i.e., $x\in \Omega\setminus\Omega_2$. The claim is proved and thus (\ref{pf-thm-Omega2}) is valid.

Now we are ready to improve the estimate for $\displaystyle \int_{\Omega} \theta_{d,m}(x) dx$. For $x\in\Omega_2$, since $\Omega_2 \subseteq \Omega_1$,  the estimate (\ref{pf-thm-Omega1-x}) still holds, i.e.,
$$
\theta_{d,m}(x) \leq m(x) + \frac{1 }{2 |\Omega|}\int_{\Omega}  \theta_{d,m}(x) dx.
$$
Then
\begin{eqnarray*}
\int_{\Omega_2} \theta_{d,m}(x) dx &\leq& \int_{\Omega} m(x) dx + \frac{|\Omega_2| }{2 |\Omega|}\int_{\Omega}  \theta_{d,m}(x) dx\\
 &\leq& \int_{\Omega} m(x) dx + \frac{1 }{2 }\int_{\Omega_2}  \theta_{d,m}(x) dx+ \frac{|\Omega_2| }{2 |\Omega|}\int_{\Omega\setminus\Omega_2}  \theta_{d,m}(x) dx,
\end{eqnarray*}
which yields that
\begin{equation}\label{pf-thm-Omega2-int}
\int_{\Omega_2} \theta_{d,m}(x) dx \leq 2\int_{\Omega} m(x) dx + \frac{|\Omega_2| }{ |\Omega|}\int_{\Omega\setminus\Omega_2}  \theta_{d,m}(x) dx.
\end{equation}
Moreover, we analyze the solution $\theta_{d,m}$ in $\Omega\setminus\Omega_2$. According to the equation (\ref{single-ss}), the estimates (\ref{pf-thm-Omega2}) and (\ref{pf-thm-Omega2-int}) and the fact that $a(x)\le 1$, one has
\begin{eqnarray*}
\int_{\Omega\setminus\Omega_2} \theta_{d,m}^2(x) dx  &=& \int_{\Omega\setminus\Omega_2}  m(x)\theta_{d,m}(x) dx + d\int_{\Omega\setminus\Omega_2}  \left( \int_{\Omega} k(x,y) \theta_{d,m}(y) dy - a(x) \theta_{d,m} (x)\right) dx\\
&\leq & K_2 d \int_{\Omega\setminus\Omega_2}  m(x) dx - d\int_{\Omega_2}  \left( \int_{\Omega} k(x,y) \theta_{d,m}(y) dy - a(x) \theta_{d,m} (x)\right) dx\\
&\leq & K_2 d\int_{\Omega} m(x) dx +  d \int_{\Omega_2} \theta_{d,m}(x) dx\\
&\leq & \left(K_2 d  + 2d + {2\over \min_{\bar\Omega} a(x) } \frac{1 }{ |\Omega|}\int_{\Omega\setminus\Omega_2}  \theta_{d,m}(x) dx\right)\int_{\Omega} m(x) dx\\
&\leq & (K_2 +2) d\int_{\Omega} m(x) dx + {2\over |\Omega|}\left({\int_{\Omega} m(x) dx\over \min_{\bar\Omega} a(x) }\right)^2 + {1\over 2}\int_{\Omega\setminus\Omega_2} \theta_{d,m}^2(x) dx .
\end{eqnarray*}
This indicates that for $d\geq 1$,
$$
\int_{\Omega\setminus\Omega_2} \theta_{d,m}^2(x) dx \leq 2 (K_2 +2) d\int_{\Omega} m(x) dx + {4 \over |\Omega|}\left({\int_{\Omega} m(x) dx\over \min_{\bar\Omega} a(x) }\right)^2 \leq K_3 d,
$$
where
$$
K_3 =  2 (K_2 +2) \int_{\Omega} m(x) dx  + {4 \over |\Omega|}\left({1\over \min_{\bar\Omega} a(x) }\right)^2 \left(\int_{\Omega} m(x) dx\right)^2.
$$
Therefore, together with (\ref{pf-thm-Omega2-int}), for $d\geq 1$
\begin{eqnarray*}
\int_{\Omega} \theta_{d,m}(x) dx &=& \int_{\Omega_2} \theta_{d,m}(x) dx +\int_{\Omega\setminus\Omega_2} \theta_{d,m}(x) dx\\
&\leq &  2\int_{\Omega} m(x) dx+ \frac{|\Omega_2| }{ |\Omega|}\int_{\Omega\setminus\Omega_2}  \theta_{d,m}(x) dx +\int_{\Omega\setminus\Omega_2} \theta_{d,m}(x) dx\\
&\leq &  2\int_{\Omega} m(x) dx + 2 \left( |\Omega\setminus\Omega_2| \int_{\Omega\setminus\Omega_2} \theta_{d,m}^2(x) dx\right)^{1\over 2}\\
&\leq &  2\left( \int_{\Omega} m(x) dx + \sqrt{K_3 |\Omega|}\right) \sqrt{d}.
\end{eqnarray*}
Set
$$
C_1 = 2\left( \int_{\Omega} m(x) dx + \sqrt{K_3 |\Omega|}\right).
$$
The desired estimate (\ref{thm-upper}) follows.
The proof is complete.
\end{proof}

\begin{remark}
It is worth pointing out that (\ref{remark-bounded-d}) indicates that the total equilibrium population is bounded if $d$ is bounded. This, together with the estimate (\ref{u-m-bound}) in the proof of Theorem \ref{thm-existence}, shows that the unboundedness of total equilibrium population under the prescribed total amount of resources is  due to the unboundedness of diffusion rate $d$ and $\|m\|_{L^{\infty}}$.
\end{remark}

\subsection{Lower bound for $\displaystyle \sup\,\left\{ \int_{\Omega}\theta_{d,m}\,dx \, \Big|\, m\in\mathcal{M}_1 \right\} $}

In this subsection, we show that there exists $C_0$, independent of $ m\in\mathcal{M}_1$,  such that
\begin{equation}\label{thm-lower}
\sup\,\left\{ \int_{\Omega}\theta_{d,m}\,dx \, \Big|\, m\in\mathcal{M}_1 \right\} \geq  C_0 \sqrt{d}.
\end{equation}

For this purpose, we construct examples for $d$ large as follows
	\begin{equation}\label{example-m-epsilon}
		m_{d}(x) =\begin{cases}
			\displaystyle 0  &x\in \Omega\setminus \Omega_{0,d},\\
			\displaystyle  M_d  &x\in \Omega_{0,d},
		\end{cases}
	\end{equation}
	where  $\Omega_{0,d} $ denotes a ball with center $x_0\in\Omega$, radius $\displaystyle(M_d \omega_n)^{-{1\over n}}$ with $M_d$  large  enough such that $\Omega_{0,d}  \subset \Omega$ and
$\omega_n$ denoting the volume of the unit ball in $\mathbb R^n$. Moreover, assume that
$$
\lim_{d\rightarrow \infty}{1\over  a(x_0)} {M_d\over d} = \alpha \in (1,\infty].
$$
Obviously $m_d \in\mathcal{M}_1$ and thus to show the lower bound (\ref{thm-lower}), it suffices to show that  there exists $C_0>0$, independent of  $d$, such that
	\begin{equation}\label{thm-lower-example}
	\int_{\Omega} \theta_{d,m_d} (x)dx \geq C_0 \sqrt{d},
	\end{equation}
where $\theta_{d,m_d}$ denotes the unique positive steady state to the problem (\ref{main single}) with  $m$ replaced by $m_d$.

\begin{proof}[Proof of the lower bound (\ref{thm-lower-example})]
First of all, it is routine to show that
\begin{equation*}
\theta_{d,m_d} (x)=\begin{cases}
\displaystyle {1\over 2} \left[ -d  a(x) + \sqrt{d^2 a^2(x) + 4 d  \int_{\Omega} k(x,y) \theta_{d,m_d} (y) dy} \right] &x\in \Omega\setminus \Omega_{0,d},\\[5mm]
\displaystyle  {1\over 2} \left[ M_d -d  a(x)  + \sqrt{\left( M_d -d a(x) \right)^2  + 4 d  \int_{\Omega} k(x,y) \theta_{d ,m_d} (y) dy} \right]   &x\in \Omega_{0,d}.
\end{cases}
\end{equation*}
Thanks to Theorem \ref{theorem-unbounded-rate}, one sees that
\begin{equation}\label{pf-theta-d}
\lim_{ d \rightarrow \infty} \frac{ \int_{\Omega} \theta_{d ,m_d}(x) dx}{ d } =0,
\end{equation}
and
\begin{equation}\label{pf-int-theta-d}
\lim_{d\rightarrow \infty} \frac{ \int_{\Omega}k(x,y) \theta_{d ,m_d}(y) dy}{ d } =0
\end{equation}
uniformly in $\Omega$.

For $\displaystyle x\in \Omega\setminus \Omega_{0,d }$, by Taylor expansion, 
we can derive
\begin{eqnarray}
\theta_{d ,m_d} (x) &=&
 {1\over 2} \left[ -d  a(x) + \sqrt{d^2 a^2(x) + 4 d  \int_{\Omega} k(x,y) \theta_{d ,m_d}(y) dy} \right]\cr
 &= &   {d \over 2}a(x)\left[ - 1 + \displaystyle\sqrt{  1 + {4\over a^2(x)}\frac{\int_{\Omega}k(x,y) \theta_{d ,m_d}(y) dy}{ d }} \right]\cr
 &=&  \frac{\displaystyle\int_{\Omega}k(x,y) \theta_{d ,m_d}(y) dy}{ a(x)} - (1+\xi)^{-{3\over 2}}a^{-3}(x){1\over d }\left(\int_{\Omega}k(x,y) \theta_{d ,m_d}(y) dy \right)^2,\nonumber
\end{eqnarray}
where
$$
0<\xi (x ) \leq {4\over a^2(x)}\frac{ \int_{\Omega}k(x,y) \theta_{d ,m_d}(y) dy}{ d }.
$$
This yields that
\begin{eqnarray}
&& \int_{\Omega\setminus \Omega_{0,d}}(1+\xi)^{-{3\over 2}}a^{-2}(x){1\over d }\left(\int_{\Omega}k(x,y) \theta_{d ,m_d}(y) dy \right)^2  dx\cr
&=&  \int_{\Omega\setminus \Omega_{0,d}}\left( \int_{\Omega}k(x,y) \theta_{d ,m_d}(y) dy\right)dx -\int_{\Omega\setminus \Omega_{0,d}} a(x) \theta_{d,m_d} (x) dx\cr
&=&   \int_{\Omega\setminus \Omega_{0,d}} \left(\int_{\Omega}k(y,x) \theta_{d,m_d}(x) dx\right) dy -\int_{\Omega\setminus \Omega_{0,d}} a(x) \theta_{d,m_d} (x) dx\cr
&= &   \int_{\Omega} \left(a(x ) - \int_{\Omega_{0,d}} k(y,x) dy \right) \theta_{d,m_d}(x) dx  -\int_{\Omega\setminus \Omega_{0,d}} a(x) \theta_{d,m_d} (x) dx \cr
&=&  \int_{\Omega_{0,d}} a(x) \theta_{d,m_d} (x) dx -  \int_{\Omega_{0,d}} \left(\int_{\Omega} k(y,x) \theta_{d,m_d}(x) dx  \right)  dy\cr
&= &  \int_{\Omega_{0,d}} {a(x) \over 2} \left[M_d -da(x) + \sqrt{\left(M_d -da(x)\right)^2  + 4 d \int_{\Omega} k(x,y) \theta_{d,m_d} (y) dy} \right] dx \cr
&& - d \int_{\Omega_{0,d}} \frac{\int_{\Omega} k(x,y) \theta_{d,m_d}(y) dy }{d}  dx.\nonumber
\end{eqnarray}
Thus, thanks to  (\ref{pf-int-theta-d}) and the assumption $\displaystyle \lim_{d\rightarrow \infty}{1\over  a(x_0)} {M_d\over d} = \alpha \in (1,\infty]$, one has
\begin{equation}\label{pf-thm-limit}
\lim_{d\rightarrow \infty} \int_{\Omega\setminus \Omega_{0,d}} (1+\xi)^{-{3\over 2}}a^{-2}(x){1\over d}\left(\int_{\Omega}k(x,y) \theta_{d,m_d}(y) dy \right)^2  dx= \left( 1-{1\over\alpha} \right)a(x_0).
\end{equation}
Notice that $\displaystyle a(x)=\int_{\Omega} k(y,x) dy$ is strictly positive and continuous  in $\bar\Omega$, and $\displaystyle \lim_{d\rightarrow +\infty}  \xi(x)=0$ uniformly in $\Omega$. Hence (\ref{pf-thm-limit}) indicates that there exists a constant $C>0$ such that for $d$ large,
\begin{equation}\label{pf-lim-(1-1/d)}
\int_{\Omega\setminus \Omega_{0,d}} {1\over d}\left(\int_{\Omega}k(x,y) \theta_{d,m_d}(y) dy \right)^2  dx \geq C.
\end{equation}

Hence under the assumption   $\displaystyle \lim_{d\rightarrow \infty}{1\over  a(x_0)} {M_d\over d} = \alpha \in (1,\infty]$, for $d$ large, one has
$$
 Cd  \leq  \int_{\Omega} \left(\int_{\Omega}k(x,y) \theta_{d ,m_d}(y) dy \right)^2 dx\leq \| k\|^2_{L^{\infty}} |\Omega| \left(\int_{\Omega}  \theta_{d ,m_d}(y) dy \right)^2.
$$
Therefore,
$$
 \int_{\Omega}  \theta_{d ,m_d}(x) dx \geq \sqrt{{C \over 2 |\Omega|} } {1\over\| k\|_{L^{\infty}}}  \sqrt{d }.
$$
and (\ref{thm-lower-example}) follows immediately by setting
$$
C_0=\sqrt{{C \over 2 |\Omega|} } {1\over\| k\|_{L^{\infty}}}.
$$
\end{proof}

The estimates (\ref{thm-upper}) and  (\ref{thm-lower-example}) yield  Theorem \ref{thm-upper-lower}.

\section{Optimal characterization of total population}
In this section, we first present the proof of Theorem \ref{thm-optimal} and then as an application of Theorem \ref{thm-optimal},  some concrete examples  are constructed  to demonstrate for $d$ large, how to choose $m$ to support more populations.

\subsection{Proof of Theorem \ref{thm-optimal}}
For  Theorem \ref{thm-optimal}, in order to prove the sufficiency of the assumption \textbf{(A)},  besides the sub-domain
    $$
    \Omega_{\varepsilon_0}=\left\{x\in\Omega \, \big | \, m(x)> d(1+\varepsilon_0)a(x) \right\}
    $$
    we also introduce
    $$
    \Omega_d = \left\{ x\in \Omega \, \big |\, m(x) \leq d^{3\over 4} a(x)  \right\}.
    $$
    Then make use of the structure of the equation satisfied by $\theta_{d,m}$, the relation between $\Omega_{\varepsilon_0}$ and $ \Omega_d$, and Taylor expansion to derive the desired conclusion. However, to  prove the necessity of the assumption \textbf{(A)}, we argue by contradiction. Suppose   there exists a sequence $d_\ell\rightarrow\infty$ as $\ell\rightarrow \infty$ and $m_{\ell}\in\mathcal{M}_1$ such that $(d_{\ell}, m_{\ell})\in \mathcal{S}$ and for any given $\varepsilon>0$,
$\displaystyle \int_{\Omega^{\ell}_{\varepsilon}} m_{\ell}(x) dx < \varepsilon
$ for $\ell$ large enough, where
$$
\Omega^{\ell}_{\varepsilon}:=\left\{x\in\Omega \, \big | \, m_{\ell}(x) > d_{\ell}(1+\varepsilon)a(x) \right\}.
$$
To derive a contradiction, the key technique is to provide a more precise estimate for $\theta_{d_{\ell},m_{\ell}}$ in $\Omega\setminus \Omega^{\ell}_{\varepsilon}$, i.e., \textbf{Claim 1} in the proof.
\begin{proof}[Proof of Theorem \ref{thm-optimal}]
First, we prove that the assumption \textbf{(A)} is  sufficient  for the total equilibrium population  reaching order $\sqrt{d}$ as $d\rightarrow \infty$.

 Based on the equation (\ref{single-ss}),
\begin{eqnarray*}
\theta_{d,m} (x) &=&
 {1\over 2} \left[m(x) -d a(x) + \sqrt{(m(x) -d a(x))^2 + 4 d \int_{\Omega} k(x,y) \theta_{d,m}(y) dy} \right]\\
&=&{da\over 2}\left[ {m\over da} - 1 + \sqrt{\left(  1-{m\over da}\right)^2 + {4\over a^2}  \frac{\int_{\Omega} k(x,y) \theta_{d,m}(y) dy}{d}} \right].
\end{eqnarray*}
Denote
$$
\Omega_d = \left\{ x\in \Omega \, \big |\, m(x) \leq d^{3\over 4} a(x)  \right\},\ \ \Omega_{\varepsilon}=\left\{x\in\Omega \, \big | \, m(x)> d(1+\varepsilon)a(x) \right\}.
$$
Obviously there exists $D_1>D$ such that for $d>D_1$, we have
$$
{m(x)\over da(x)}<{1\over 2}\ \ \textrm{for}\  x\in\Omega_d.
$$
Then by Taylor expansion $\sqrt{\beta^2+z}=\beta+\frac{1}{2}\beta^{-1}z-\frac{1}{8}(\beta^2+\xi)^{-\frac{3}{2}}z^2,~\beta>0$, it is standard to derive that for $x\in\Omega_d$, $d>D_1$,
\begin{eqnarray*}
\theta_{d,m}(x) &= &  \left(a- {m\over d}\right)^{-1}\int_{\Omega} k(x,y) \theta_{d,m}(y) dy -\left(\left(1- {m\over da}\right)^{2}+\xi\right)^{-{3\over 2}}{1\over a^3 d}\left(\int_{\Omega} k(x,y) \theta_{d,m}(y) dy \right)^2,
\end{eqnarray*}
where
$$
0<\xi < {4\over a^2}  \frac{\int_{\Omega} k(x,y) \theta_{d,m}(y) dy}{d}.
$$
This implies that
\begin{eqnarray}\label{pf-thm-optimal-1-taylor}
&&\int_{\Omega_d}\left(a- {m\over d}\right)\left(\left(1- {m\over da}\right)^{2}+\xi\right)^{-{3\over 2}}{1\over a^3d}\left(\int_{\Omega} k(x,y) \theta_{d,m}(y) dy \right)^2dx \cr
&=& \int_{\Omega_d} \left({m(x)\over d}\theta_{d,m}(x) +\int_{\Omega} k(x,y) \theta_{d,m}(y) dy -a (x)\theta_{d,m} (x) \right)dx
\end{eqnarray}

Let us estimate the right hand side of (\ref{pf-thm-optimal-1-taylor}) first. Notice that
for $d\geq 1$, $\Omega_{\varepsilon_0}\subseteq \Omega\setminus\Omega_{d}$, and $$|\Omega\setminus\Omega_{d}|\le \int_{\{m\ge d^{\frac{3}{4}}a(x)\}}\frac{m(x)}{d^{\frac{3}{4}}\min_{\bar{\Omega}} a} dx\leq \frac{1}{\min_{\bar\Omega} a}d^{-{3\over 4}},$$
since $m\in\mathcal{M}_1$.
Then it follows from Theorem \ref{theorem-unbounded-rate} and the assumption \textbf{(A)} that for $d>D_1$
\begin{eqnarray}
&& \int_{\Omega_d} \left({m(x)\over d}\theta_{d,m}(x) +\int_{\Omega} k(x,y) \theta_{d,m}(y) dy -a(x) \theta_{d,m} (x) \right)dx\cr
&\geq & \int_{\Omega_d} \left( \int_{\Omega} k(x,y) \theta_{d,m}(y) dy -a(x) \theta_{d,m} (x) \right)dx\cr
&=&   \int_{\Omega\setminus\Omega_d} \left(  a(x) \theta_{d,m}(x) -\int_{\Omega} k(x,y) \theta_{d,m}(y) dy \right)dx\cr
&\geq &\int_{\Omega_{\varepsilon_0}} a(x) \theta_{d,m}(x)dx - \int_{\Omega\setminus\Omega_d}\int_{\Omega} k(x,y) \theta_{d,m}(y) dy dx\cr
&=& \int_{\Omega_{\varepsilon_0}} {a\over 2} \left[m -d a + \sqrt{(m -d a)^2 + 4 d \int_{\Omega} k(x,y) \theta_{d,m}(y) dy} \right] dx- \int_{\Omega\setminus\Omega_d}\int_{\Omega} k(x,y) \theta_{d,m}(y) dy dx\cr
&\geq &\int_{\Omega_{\varepsilon_0}} a(m-da)dx - C_1\|k \|_{L^{\infty}}\sqrt{d} |\Omega\setminus\Omega_d|\cr
&\geq &\int_{\Omega_{\varepsilon_0}} a\left(m- {m\over 1+\varepsilon_0}\right)dx - C_1\|k \|_{L^{\infty}} \frac{1}{\min_{\bar\Omega} a}d^{-{1\over 4}}\cr
&\geq &{\varepsilon_0^2 \over 1+\varepsilon_0}\min_{\bar\Omega} a - C_1\|k \|_{L^{\infty}} \frac{1}{\min_{\bar\Omega} a}d^{-{1\over 4}},\nonumber
\end{eqnarray}
This, together with (\ref{pf-thm-optimal-1-taylor}), indicates that there exists $D_2>D_1$ such that
\begin{eqnarray*}
{1\over 2}{\varepsilon_0^2 \over 1+\varepsilon_0}\min_{\bar\Omega} a &\leq &   \int_{\Omega_d} \left({m(x)\over d}\theta_{d,m}(x) +\int_{\Omega} k(x,y) \theta_{d,m}(y) dy -a(x) \theta_{d,m} (x) \right)dx\\
&= &\int_{\Omega_d}\left(a- {m\over d}\right)\left(\left(1- {m\over da}\right)^{2}+\xi\right)^{-{3\over 2}}{1\over a^3  d}\left(\int_{\Omega} k(x,y) \theta_{d,m}(y) dy \right)^2dx\\
&\le &\int_{\Omega_d}\left(\frac{1}{4}\right)^{-{3\over 2}}{1\over a^2d} \left(\int_{\Omega} k(x,y) \theta_{d,m}(y) dy \right)^2dx\\
&\leq &8 \left(\min_{\bar\Omega} a \right)^{-2}\|k \|^2_{L^{\infty}} |\Omega| {1\over d} \left(\int_{\Omega}\theta_{d,m}\,dx\right)^2.
\end{eqnarray*}
It is routine to check that the lower bound estimate$\displaystyle \int_{\Omega}\theta_{d,m}\,dx \geq  C \sqrt{d}$ is valid with
$$
C=\frac{\varepsilon_0\left(\min_{\bar\Omega} a \right)^{\frac{3}{2}}}{4\|k \|_{L^{\infty}}\sqrt{(1+\varepsilon_0)|\Omega|}}.
$$

Next, we prove the  necessity of the assumption \textbf{(A)} for the total equilibrium population  reaching order $\sqrt{d}$ as $d\rightarrow \infty$ in the subset $\mathcal{S}\subseteq  \{d>0 \}\times \mathcal{M}_1$.
Suppose the assumption \textbf{(A)} does not hold in $\mathcal{S}$, i.e., there exists a sequence $d_\ell\rightarrow\infty$ as $\ell\rightarrow \infty$ and $m_{\ell}\in\mathcal{M}_1$ such that $(d_{\ell}, m_{\ell})\in \mathcal{S}$ and for any given $\varepsilon>0$,
$\displaystyle \int_{\Omega^{\ell}_{\varepsilon}} m_{\ell}(x) dx < \varepsilon
$ for $\ell$ large enough, where
$$
\Omega^{\ell}_{\varepsilon}:=\left\{x\in\Omega \, \big | \, m_{\ell}(x) > d_{\ell}(1+\varepsilon)a(x) \right\}.
$$
It suffices to show
\begin{equation}\label{pf-thm-0}
\lim_{\ell \rightarrow \infty} \frac{\int_{\Omega}\theta_{d_{\ell},m_{\ell}}\,dx}{\sqrt{d_{\ell}}}=0.
\end{equation}
For this purpose, the following claim is crucial:

\textbf{Claim 1.} {\it There exists $L_1>0$, such that for $\ell>L_1$,  $\displaystyle\theta_{d_{\ell},m_{\ell}}(x) < 2\left( \max_{\bar\Omega} a\right) d_{\ell}\varepsilon$ in $\Omega\setminus \Omega^{\ell}_{\varepsilon}$.}

Assume that the claim is true.  Based on the equation satisfied by $\theta_{d_{\ell}}$, one has for $\ell>L_1$
\begin{eqnarray}\label{pf-thm-theta2}
&& \int_{\Omega\setminus \Omega^{\ell}_{\varepsilon}}\theta_{d_{\ell},m_{\ell}}^2 dx\cr
 &=&\int_{\Omega\setminus \Omega^{\ell}_{\varepsilon}} m_{\ell} \theta_{d_{\ell},m_{\ell}} dx +d_{\ell}\int_{\Omega\setminus \Omega^{\ell}_{\varepsilon}}\left(\int_{\Omega} k(x,y) \theta_{d_{\ell},m_{\ell}}(y) dy -a(x) \theta_{d_{\ell},m_{\ell}} (x)\right)dx\cr
&=& \int_{\Omega\setminus \Omega^{\ell}_{\varepsilon}} m_{\ell} \theta_{d_{\ell},m_{\ell}} dx - d_{\ell}\int_{\Omega^{\ell}_{\varepsilon}}\left(\int_{\Omega} k(x,y) \theta_{d_{\ell},m_{\ell}}(y) dy -a(x) \theta_{d_{\ell},m_{\ell}} (x)\right)dx\cr
&\leq &\int_{\Omega\setminus \Omega^{\ell}_{\varepsilon}} m_{\ell} \theta_{d_{\ell},m_{\ell}} dx + d_{\ell}\int_{\Omega^{\ell}_{\varepsilon}}a(x) \theta_{d_{\ell},m_{\ell}} (x) dx\cr
&\leq & 2\left( \max_{\bar\Omega} a\right) d_{\ell}\varepsilon + d_{\ell}\int_{\Omega^{\ell}_{\varepsilon}}{a  \over 2} \left(m_{\ell}(x) -d_{\ell} a + \sqrt{(m _{\ell}-d_{\ell} a)^2 + 4 d_{\ell} \int_{\Omega} k(x,y) \theta_{d_{\ell},m_{\ell}}(y) dy} \right) dx\cr
&\leq & 2\left( \max_{\bar\Omega} a\right) d_{\ell}\varepsilon + d_{\ell}\int_{\Omega^{\ell}_{\varepsilon}}{a  \over 2} \left(m_{\ell}(x) -d_{\ell} a +\sqrt{(m_{\ell} -d_{\ell} a)^2 }+ \sqrt{ 4 d_{\ell} \int_{\Omega} k(x,y) \theta_{d_{\ell},m_{\ell}}(y) dy} \right) dx\cr
&\leq &2\left( \max_{\bar\Omega} a\right) d_{\ell}\varepsilon + d_{\ell}\int_{\Omega^{\ell}_{\varepsilon}} a(x) \left(m_{\ell}(x)  + \sqrt{ d_{\ell} \int_{\Omega} k(x,y) \theta_{d_{\ell},m_{\ell}}(y) dy} \right) dx\cr
&\leq &2\left( \max_{\bar\Omega} a\right) d_{\ell}\varepsilon +\left( \max_{\bar\Omega} a\right) d_{\ell}\varepsilon+
 \left(C_1 \|k\|_{L^{\infty}}\right)^{1\over 2} \left( \max_{\bar\Omega} a\right)d_{\ell}^{7\over 4} |\Omega^{\ell}_{\varepsilon} |
 \cr
&\leq &3\left( \max_{\bar\Omega} a\right) d_{\ell}\varepsilon + 
 \left(C_1 \|k\|_{L^{\infty}}\right)^{1\over 2} \left( \max_{\bar\Omega} a\right)d_{\ell}^{7\over 4} \displaystyle\int_{\Omega^{\ell}_{\varepsilon}}\frac{m_{\ell}(x)}{(1+\varepsilon)d_\ell a(x)} dx
  \cr
&\leq &3\left( \max_{\bar\Omega} a\right) d_{\ell}\varepsilon + 
 \left(C_1 \|k\|_{L^{\infty}}\right)^{1\over 2} \left( \max_{\bar\Omega} a\right)d_{\ell} {\varepsilon\over \min_{\bar\Omega} a}:=C(a,k) d_{\ell} \varepsilon
\end{eqnarray}
where the fourth last inequality  is due to Theorem  \ref{theorem-unbounded-rate}. Moreover according to the definition of $\Omega^{\ell}_{\varepsilon}$, $m\in\mathcal{M}_1$ and Theorem  \ref{theorem-unbounded-rate}, it is easy to check that
\begin{eqnarray*}
 \int_{\Omega^{\ell}_{\varepsilon}} \theta_{d_{\ell},m_{\ell}} (x) dx &=&\int_{\Omega^{\ell}_{\varepsilon}}{1  \over 2} \left(m_{\ell}(x) -d_{\ell} a + \sqrt{(m_{\ell} -d_{\ell} a)^2 + 4 d_{\ell} \int_{\Omega} k(x,y) \theta_{d_{\ell},m_{\ell}}(y) dy} \right) dx\\
&\leq & \int_{\Omega^{\ell}_{\varepsilon}} \left(m_{\ell}(x)  + \sqrt{ d_{\ell} \int_{\Omega} k(x,y) \theta_{d_{\ell},m_{\ell}}(y) dy} \right) dx\\
& \leq&  \varepsilon+ \left(C_1 \|k\|_{L^{\infty}}\right)^{1\over 2} \left(\min_{\bar\Omega} a\right)^{-1}d_{\ell}^{-{1\over 4}}
\end{eqnarray*}
This, together with (\ref{pf-thm-theta2}), implies that there exists $L_2 >L_1$ such that for $\ell >L_2$
\begin{eqnarray*}
 \int_{\Omega} \theta_{d_{\ell},m_{\ell}} (x) dx &=& \int_{\Omega^{\ell}_{\varepsilon}} \theta_{d_{\ell},m_{\ell}} (x) dx+ \int_{\Omega\setminus \Omega^{\ell}_{\varepsilon}}\theta_{d_{\ell},m_{\ell}} (x)dx\\
&\leq&  \varepsilon+ \left(C_1 \|k\|_{L^{\infty}}\right)^{1\over 2} \left(\min_{\bar\Omega} a\right)^{-1}d_{\ell}^{-{1\over 4}}+ |\Omega|^{1\over 2}\left(\int_{\Omega\setminus \Omega^{\ell}_{\varepsilon}}\theta_{d_{\ell},m_{\ell}}^2 dx\right)^{1\over 2}\\
&\leq & 1+ |\Omega|^{1\over 2}\left(C(a,k) d_{\ell} \varepsilon\right)^{1\over 2}.
\end{eqnarray*}
By arbitrariness of $\varepsilon$, we have
$$
\lim_{d_{\ell} \rightarrow \infty} \frac{\int_{\Omega}\theta_{d_{\ell},m_{\ell}}\,dx}{\sqrt{d_{\ell}}}=0.
$$
This verifies (\ref{pf-thm-0}) and the necessity of the assumption \textbf{(A)}.

Now it remains  to verify {\bf Claim 1}. It is equivalent to show that there exists $L_1>0$,  for  $\ell>L_1$, $\displaystyle\theta_{d_{\ell}}(x) \geq  2\left( \max_{\bar\Omega} a\right) d_{\ell}\varepsilon$, $x\in\Omega$ implies that $m_{\ell}(x)>  d_{\ell} a(x) \left(1+ \varepsilon\right).
$

Set   $\displaystyle\alpha :=2\left(\max_{\bar\Omega}a \right)\varepsilon$ for clarity.
Notice that if $\theta_{d_{\ell}}(x) \geq d_{\ell}\alpha$ and $m_{\ell}(x)< d_{\ell} a(x)$ happen together, then due to the equation satisfied by $\theta_{d_{\ell}}$ and Theorem \ref{theorem-unbounded-rate}, one has
\begin{eqnarray*}
 d_{\ell}\alpha \leq \theta_{d_{\ell}}(x)&=&{1\over 2} \left[m_{\ell}(x) -d_{\ell} a(x) + \sqrt{(m_{\ell}(x) -d_{\ell} a(x))^2 + 4 d_{\ell} \int_{\Omega} k(x,y) \theta_{d_{\ell},m_{\ell}}(y) dy} \right]\\
&=& \frac{  2d_{\ell} \int_{\Omega} k(x,y) \theta_{d_{\ell},m_{\ell}}(y) dy}{  -m_{\ell}(x) +d_{\ell} a(x) + \sqrt{(m_{\ell}(x) -d_{\ell} a(x))^2 + 4 d_{\ell} \int_{\Omega} k(x,y) \theta_{d_{\ell},m_{\ell}}(y) dy}}\\
&\leq &\sqrt{ d_{\ell} \int_{\Omega} k(x,y) \theta_{d_{\ell},m_{\ell}}(y) dy}\\
&\leq &\left(C_1 \|k\|_{L^{\infty}}\right)^{1\over 2}d_{\ell}^{3\over 4}.
\end{eqnarray*}
This yields that
$$
d_{\ell}< \frac{\left(C_1 \|k\|_{L^{\infty}}\right)^2}{\alpha^4}.$$
Hence, if $\displaystyle d_{\ell}\ge \frac{\left(C_1 \|k\|_{L^{\infty}}\right)^2}{\alpha^4}$, then $\theta_{d_{\ell}}(x) \geq d_{\ell}\alpha$ guarantees that
\begin{equation}\label{pf-thm-estimate}
m_{\ell}(x)\geq d_{\ell} a(x).
\end{equation}

Moreover,  if $\theta_{d_{\ell}}(x) \geq d_{\ell}\alpha$, direct calculation gives
\begin{eqnarray*}
 d_{\ell}\alpha \leq \theta_{d_{\ell}}(x) &=&{1\over 2} \left[m_{\ell}(x) -d_{\ell} a(x) + \sqrt{(m_{\ell}(x) -d_{\ell} a(x))^2 + 4 d_{\ell} \int_{\Omega} k(x,y) \theta_{d_{\ell},m_{\ell}}(y) dy} \right]\\
&\leq &    \sqrt{(m_{\ell}(x) -d_{\ell} a(x))^2 + 4 d_{\ell} \int_{\Omega} k(x,y) \theta_{d_{\ell},m_{\ell}}(y) dy},
\end{eqnarray*}
which, together with Theorem \ref{theorem-unbounded-rate}, yields that
\begin{equation}\label{pf-thm-estimate-2}
\left({m_\ell(x)\over d_{\ell}} - a(x)\right)^2\geq \alpha^2 - \frac{4\int_{\Omega} k(x,y) \theta_{d_{\ell},m_{\ell}}(y) dy}{d_{\ell}}\geq \alpha^2 - \frac{4C_1\|k\|_{L^{\infty}}}{\sqrt{d_{\ell}}}>{\alpha^2\over 4},
\end{equation}
provided that $d_\ell$ is  sufficiently large. At the end, choose  $L_1>0$ large  enough  such that for $\ell>L_1$, $\displaystyle d_{\ell}\ge \frac{\left(C_1 \|k\|_{L^{\infty}}\right)^2}{\alpha^4}$ and (\ref{pf-thm-estimate-2}) holds. Then by (\ref{pf-thm-estimate}) and (\ref{pf-thm-estimate-2}), we can further obtain
$$
{m_{\ell}(x)\over d_{\ell}} > a(x)+{\alpha\over 2}\geq a(x)\left(1+ \frac{\alpha}{ 2\max_{\bar\Omega}a} \right)=a(x)\left(1+\varepsilon\right).
$$
The proof is complete.
\end{proof}

\subsection{Applications}
Based on the equivalent   criterion  established in Theorem \ref{thm-optimal},  a series examples are constructed to show that how the concentration of resources, including the height and locations, and properties of nonlocal kernel functions affect the total population for large diffusion rate.

\medskip

\textbf{Example 1.} Define
	\begin{equation*}
		m_{d}^{\alpha,\beta}(x) =\begin{cases}
			\displaystyle 0  &x\in \Omega\setminus \Omega^{\alpha, \beta }_{d},\\
			\displaystyle  \alpha d^{\beta}   &x\in \Omega^{\alpha,\beta}_{d},
		\end{cases}
	\end{equation*}
	where  $\Omega^{\alpha,\beta}_{d} \subseteq \Omega$ with $\displaystyle |\Omega^{\alpha,\beta}_{d}|=  (\alpha d^{\beta})^{-1},$ $\alpha, \beta>0$. Obviously $m^{\alpha,\beta}_d \in\mathcal{M}_1$. Let $\theta_{d}^{\alpha,\beta}$ denote the unique positive steady state to the problem (\ref{main single}) with  $m(x)$ replaced by $m_{d}^{\alpha,\beta}(x)$. Thanks to Theorem \ref{thm-optimal}, we have the following statements.
\begin{itemize}
\item[(i)] If $0<\beta<1$, then  for any $\alpha>0$, we have
$\displaystyle
\lim_{d\rightarrow \infty}\frac{\int_{\Omega}\theta_{d}^{\alpha,\beta} dx }{\sqrt{d}} =0.
$
\item[(ii)] If $\beta>1$, then for any $\alpha>0$, we have $\displaystyle \int_{\Omega}\theta_{d}^{\alpha,\beta} dx$ is of order $\sqrt d$ as $d\rightarrow \infty$.
\item[(iii)] For the critical case $\beta =1$, recall that
 $$
 0<\min_{\bar\Omega} a(x)\leq \max_{\bar\Omega} a(x)\leq 1, \ \textrm{where}\ a(x) = \int_{\Omega} k(y,x) dy.
 $$
 Then
\begin{itemize}
\item  $\displaystyle
\lim_{d\rightarrow \infty}\frac{\int_{\Omega}\theta_{d}^{\alpha,1} dx }{\sqrt{d}} =0
$
when $\displaystyle 0<\alpha\leq \min_{\bar\Omega} a(x)$,
\item
   $\displaystyle \int_{\Omega}\theta_{d}^{\alpha,1} dx$ is of order $\sqrt d$ as $d\rightarrow \infty$ when $\displaystyle \alpha>\max_{\bar\Omega} a(x)$.
\end{itemize}
\end{itemize}
Notice that for the cases discussed  in \textbf{Example 1}, we only require that $\Omega^{\alpha,\beta}_{d}$ is a measurable subset in $\Omega$ with $\displaystyle |\Omega^{\alpha,\beta}_{d}|=  (\alpha d^{\beta})^{-1}.$

\medskip

If in addition, assume that $\displaystyle\min_{\bar\Omega} a(x) < \max_{\bar\Omega} a(x)$, then the case that
$$
 \min_{\bar\Omega} a(x) < \alpha \leq \max_{\bar\Omega} a(x)\leq 1, \ \beta=1
 $$
 is not mentioned in \textbf{Example 1}.
 Indeed, in this case, the locations where the resources concentrate, i.e., $\Omega^{\alpha,1}_{d}$, and where $\alpha>a(x)$ will affect the order of total population as $d\rightarrow \infty$.
 To better elaborate  this point, we construct an example  under some extra assumptions.

\medskip

\textbf{Example 2.} Assume that $\Omega$ is open and convex in $\mathbb R^n$, the function $J$ satisfies
\medskip\\
\noindent\textbf{(J)} \   $J(z)\in C(\mathbb R^n)$ is nonnegative, radially symmetric,  $J(0)>0$ and  $\int_{\mathbb R^n} J(z)dz =1$.
\medskip\\
Also, assume that  $J$ is compactly supported and $\textrm{diam}\, \{J>0\} \ll 1$.

Define
\begin{equation*}
		m_{d,x_0}(x) =\begin{cases}
			\displaystyle 0  &x\in \Omega\setminus \Omega_{d,x_0},\\
			\displaystyle   d  &x\in \Omega_{d,x_0},
		\end{cases}
	\end{equation*}
where $x_0\in\Omega$, $\Omega_{d,x_0}$ denotes a ball centered at $x_0$ with $|\Omega_{d,x_0}|= d^{-1}$. Let $\theta_{d,x_0}$ denote the unique positive steady state to the problem (\ref{main single}) with  $m(x)$ replaced by $m_{d,x_0}(x)$ and the kernel function $k(x,y)$ replaced by $J(x-y)$. Then
\begin{itemize}
\item
   $\displaystyle \int_{\Omega}\theta_{d,x_0} dx$ is of order $\sqrt d$ as $d\rightarrow \infty$ when $x_0$ is close to the boundary of $\Omega$,
\item  $\displaystyle
\lim_{d\rightarrow \infty}\frac{\int_{\Omega}\theta_{d,x_0} dx }{\sqrt{d}} =0
$
when $x_0$ is away from the boundary of $\Omega$.
\end{itemize}
This follows immediately from Theorem \ref{thm-optimal} and the observation that
$a(x)<1$ if $\textrm{dist} \{x,\partial\Omega\}<\textrm{diam}\, \{J>0\} /4$ and
$a(x)=1$  if $\textrm{dist} \{x,\partial\Omega\}>\textrm{diam}\, \{J>0\} /2$.

\medskip

Contrary to \textbf{Example 2}, to support more population for $d$ large, under certain assumptions,  the resources need concentrate away from the boundary. An example is constructed as follows.

\medskip

\textbf{Example 3.} Assume that $\Omega:=B_1(0)$, the function $J$ satisfies the assumption  \textbf{(J)} and  for some small $\delta>0$, we assume in addition that  $J(z)\equiv \delta, z\in B_1(0)$, $J(z)$ is strictly increasing in $|z|$ for $1<|z|\le2$.
It easy to check that $a(x)$ is radially symmetric and strictly increasing in $|x|$ for $|x|\le 1$.

Fix $a(0)<\hat\alpha<a(x)$ with $|x|=1$ and
define
	\begin{equation*}
		m_{d,x_0}^{\hat\alpha}(x) =\begin{cases}
			\displaystyle 0  &x\in \Omega\setminus \Omega^{\hat\alpha }_{d,x_0},\\
			\displaystyle  \hat\alpha d   &x\in \Omega^{\hat\alpha}_{d,x_0},
		\end{cases}
	\end{equation*}
	where $x_0\in\Omega$, $\Omega^{\hat\alpha}_{d,x_0}$ denotes a ball centered at $x_0$ with $\displaystyle |\Omega^{\hat\alpha}_{d,x_0}|=  (\hat\alpha d)^{-1}.$ Let $\theta^{\hat\alpha}_{d,x_0}$ denote the unique positive steady state to the problem (\ref{main single}) with  $m(x)$ replaced by $m^{\hat\alpha}_{d,x_0}(x)$ and the kernel function $k(x,y)$ replaced by $J(x-y)$. Then thanks to Theorem \ref{thm-optimal}, it follows that
\begin{itemize}
\item
   $\displaystyle \int_{\Omega}\theta^{\hat\alpha}_{d,x_0} dx$ is of order $\sqrt d$ as $d\rightarrow \infty$ when $x_0$ is close to the origin $0$,
\item  $\displaystyle
\lim_{d\rightarrow \infty}\frac{\int_{\Omega}\theta^{\hat\alpha}_{d,x_0} dx }{\sqrt{d}} =0
$
when $x_0$ is close to the boundary of $\Omega$.
\end{itemize}

\end{document}